\documentclass[12pt,reqno]{article}

\usepackage{apacite}
\usepackage[usenames]{color}
\usepackage{amssymb}
\usepackage{amsmath}
\usepackage{amsthm}
\usepackage{amsfonts}
\usepackage{amscd}
\usepackage{graphicx}

\usepackage{color}
\usepackage{fullpage}
\usepackage{float}

\usepackage{cases}

\usepackage{graphics,amsmath,amssymb}
\usepackage{amsthm}
\usepackage{amsfonts}
\usepackage{latexsym}

\begin{document}

\theoremstyle{plain}
\newtheorem{theorem}{Theorem}
\newtheorem{corollary}[theorem]{Corollary}
\newtheorem{lemma}[theorem]{Lemma}
\newtheorem{proposition}[theorem]{Proposition}

\theoremstyle{definition}
\newtheorem{definition}[theorem]{Definition}
\newtheorem{example}[theorem]{Example}
\newtheorem{conjecture}[theorem]{Conjecture}

\theoremstyle{remark}
\newtheorem{remark}[theorem]{Remark}

\begin{center}
\vskip 1cm{\LARGE\bf An algorithm for the Faulhaber polynomials}
\vskip 1cm \large Jos\'{e} Luis Cereceda \\
Collado Villalba, 28400 Madrid, Spain \\
{\tt jl.cereceda@movistar.es}
\end{center}

\vskip .2 in

\begin{abstract}
Let $S_p(n)$ denote the sum of $p$th powers of the first $n$ positive integers $1^p + 2^p + \cdots + n^p$. In this paper, first we express $S_p(n)$ in the so-called Faulhaber form, namely, as an even or odd polynomial in $(n + 1/2)$, according as $p$ is odd or even. Then, using the relation $S_p(n) - S_p(n-1) = n^p$, we derive a recursive formula for the associated Faulhaber coefficients. Applying Cramer's rule to the corresponding system of equations, we obtain an explicit determinant formula for the said coefficients. Furthermore, we show how to convert the (even or odd) Faulhaber polynomials in $(n+ 1/2)$ into polynomials in $S_1(n)$ for any arbitrary $p$, and vice versa.
\end{abstract}

\section{Introduction}

For integers $p \geq 0$ and $n \geq 1$, let $S_p(n)$ denote the sum of $p$th powers of the first $n$ positive integers
\begin{equation*}
S_p(n) = 1^p + 2^p + \cdots + n^p,
\end{equation*}
where it is understood that $S_p(0) =0$ for all $p$. As is well known, $S_p(n)$ can be expressed as a polynomial in $n$ of degree $p+1$ without constant term according to the Bernoulli formula (see, e.g., \citeA{sherwood} and \citeA{wu})
\begin{equation*}
S_p(n) = \frac{1}{p+1} \sum_{j=1}^{p+1} \binom{p+1}{j} (-1)^{p+1-j} B_{p+1-j} n^j, \quad p \geq 0,
\end{equation*}
where $B_0 =1$, $B_1 = -1/2$, $B_3 =0$, $B_4 = -1/30$, etc., are the Bernoulli numbers \cite{apostol}. It is relatively less well known that, for even $p=2k$, $k \geq 1$, $S_{2k}(n)$ admits the polynomial representation
\begin{equation}\label{even}
S_{2k}(n) = \sum_{m=0}^{k} f^{(2k)}_{m} \left( n+\frac{1}{2} \right)^{2m+1},
\end{equation}
while, for odd $p = 2k+1$, $k \geq 0$, $S_{2k+1}(n)$ takes the form
\begin{equation}\label{odd}
S_{2k+1}(n) = c_{2k+1} + \sum_{m=0}^{k} f^{(2k+1)}_{m} \left( n+\frac{1}{2} \right)^{2m+2},
\end{equation}
so that $S_{2k}(n)$ [$S_{2k+1}(n)$] can be expressed as an odd [even] polynomial in $n + \frac{1}{2}$. Clearly, since $S_{2k+1}(0) =0$, from \eqref{odd} we quickly obtain
\begin{equation}\label{ind}
c_{2k+1} = -\sum_{m=0}^{k}\frac{f_{m}^{(2k+1)}}{4^{m+1}}.
\end{equation}
We will refer to the polynomial forms in \eqref{even} and \eqref{odd} as Faulhaber polynomials, after the German mathematician Johann Faulhaber (1580-1635) who was the first to express $S_p(n)$ equivalently as polynomials in $S_1(n)$ for general $p$ \cite{edwards1,edwards2,knuth,beardon,dubeau,cere}.

\citeA{kelly} described a well-known method for finding $S_p(n)$ in terms of the earlier power sums $S_0(n), S_1(n), \ldots S_{p-1}(n)$ (see also \citeA{acu} for a refinement of this method showing how to obtain $S_p(n)$ from $S_{p-2}(n), S_{p-4}(n), \, \ldots\,\,$). On the other hand, there exists a procedure which allows one to deduce the coefficients $a_j^{(p)}$, $j=1,2,\ldots,p+1$, of $S_p(n) = a_{p+1}^{(p)} n^{p+1} + a_{p}^{(p)} n^{p} + \cdots + a_{1}^{(p)} n$ once the coefficients $a_i^{(p-1)}$, $i=1,2,\ldots,p$, of $S_{p-1}(n) = a_{p}^{(p-1)} n^{p} + a_{p-1}^{(p-1)} n^{p-1} + \cdots + a_{1}^{(p-1)} n$ are known; see, e.g., \citeA{budin}, \citeA{carchidi}, \citeA{owens}, \citeA{bloom}, and \citeA{torabi}. Moreover, \citeA{schultz} used the relation $S_k(n) - S_k(n-1) = n^k$ to generate recursively the coefficients $A_{k+1}, A_{k}, \ldots, A_{1}$ of $S_k(n) = A_{k+1}n^{k+1} + A_k n^k + \cdots + A_ 1 n$ for any arbitrary $k$. Note that this latter method enables us to get $S_k(n)$ without having to know the preceding sums $S_p(n)$, $p < k$ (see also \citeA{scott}, \citeA{burrows}, \citeA{bruyn}, \citeA{molnar}, and the recent paper by \citeA{marko}). Lastly, we mention that a review of these and other techniques for computing the coefficients of the power sum polynomials may be found in \citeA{kotiah} and \citeA{tanton}.

In this paper, we take up Schultz's approach to derive a recursive formula for the Faulhaber coefficients $\{ f^{(2k)}_{m}\}_{m=0}^{k}$ and $\{f^{(2k+1)}_{m} \}_{m=0}^{k}$. Then we apply Cramer's rule to the corresponding system of equations to obtain a determinant formula for each $f^{(2k)}_{m}$ and $f^{(2k+1)}_{m}$. Specifically, as we shall see, both $f^{(2k)}_{k-j}$ and $f^{(2k+1)}_{k-j}$, $j =0,1,\ldots,k$, can be expressed as a numerical factor depending on $k$ and $j$ times a determinant of order $j$ involving binomial coefficients. As an example illustrating the method presented here, we determine the Faulhaber polynomials \eqref{even} and \eqref{odd} for the case $k =5$. Furthermore, in the last part of the paper we show how the said Faulhaber polynomials \eqref{even} and \eqref{odd} can be transformed into polynomials in $S_1(n) = \frac{1}{2}n(n+1)$ for any arbitrary $k \geq 1$, and vice versa.

\section{Determination of the Faulhaber coefficients}

In order to derive a recursive formula for the coefficients $f^{(2k)}_{k}, f^{(2k)}_{k-1}, \ldots, f^{(2k)}_{0}$, we start from the relation
\begin{equation*}
S_{2k}(n) - S_{2k}(n-1) = n^{2k}.
\end{equation*}
Thus, from \eqref{even}, we have
\begin{equation}\label{even1}
\sum_{m=0}^{k}  f^{(2k)}_{m} \left[\left( n+\frac{1}{2} \right)^{2m+1} - \left( n-\frac{1}{2} \right)^{2m+1} \right] = n^{2k}.
\end{equation}
Applying the binomial theorem, the left-hand side of \eqref{even1} can be written as
\begin{equation*}
\sum_{m=0}^{k}  f^{(2k)}_{m} \sum_{j=0}^{2m+1} \frac{2^{j}}{2^{2m+1}} \binom{2m+1}{j} \big( 1- (-1)^{2m+1-j} \big) n^j.
\end{equation*}
Hence, noting that
\begin{equation*}
1- (-1)^{2m+1-j} = \begin{cases}
0,  &\text{for odd  $j$,} \\
2,  &\text{for even $j$,}
\end{cases}
\end{equation*}
it follows that
\begin{equation*}
\sum_{m=0}^{k}  f^{(2k)}_{m} \sum_{j=0}^{m} \frac{2^{2j}}{2^{2m}} \binom{2m+1}{2j} n^{2j} =
\sum_{j=0}^{k} \left[ \sum_{m = j}^{k} 4^{j-m} \binom{2m+1}{2j} f^{(2k)}_{m} \right]  n^{2j} = n^{2k}.
\end{equation*}
Equating like terms on both sides of the rightmost equality yields the following linear triangular system of $k+1$ equations in the unknowns $f^{(2k)}_{k}, f^{(2k)}_{k-1}, \ldots, f^{(2k)}_{0}$
\begin{equation}\label{recur1}
\sum_{m = j}^{k} 4^{j-m} \binom{2m+1}{2j} f^{(2k)}_{m} = \delta_{j,k}, \quad j =0,1,\ldots,k,
\end{equation}
where $\delta_{j,k}$ is a Kronecker's delta. This system can be written in matrix form as
\begin{equation*}
\setlength{\arraycolsep}{3.pt}
\begin{pmatrix}
\binom{2k+1}{1} & 0 & 0 & \hdots & 0 & 0 \\[6pt]
4^{-1}\binom{2k+1}{3} & \binom{2k-1}{1} & 0 & \hdots & 0 & 0 \\[8pt]
4^{-2}\binom{2k+1}{5} & 4^{-1}\binom{2k-1}{3} & \binom{2k-3}{1} & \hdots  & 0 & 0 \\[6pt]
\vdots & \vdots & \vdots & \ddots & \vdots & \vdots \\[7pt]
4^{-k+1}\binom{2k+1}{2k-1} & 4^{-k+2}\binom{2k-1}{2k-3} & 4^{-k+3}\binom{2k-3}{2k-5} & \hdots & \binom{3}{1} & 0 \\[7pt]
4^{-k}\binom{2k+1}{2k+1} & 4^{-k+1}\binom{2k-1}{2k-1} & 4^{-k+2}\binom{2k-3}{2k-3} & \hdots & 4^{-1}\binom{3}{3}
& \binom{1}{1} \\[5pt]
\end{pmatrix}
\begin{pmatrix} f_{k}^{(2k)} \\[6pt]  f_{k-1}^{(2k)} \\[6pt]  f_{k-2}^{(2k)} \\[2pt] \vdots \\[6pt]
  f_{1}^{(2k)} \\[6pt] f_{0}^{(2k)}
\end{pmatrix}
= \begin{pmatrix} 1\\[7pt] 0 \\[7pt] 0 \\[7pt] \vdots \\[7pt] 0 \\[7pt] 0
\end{pmatrix}.
\end{equation*}
Let us call the above $(k+1)\times (k+1)$ matrix $M_k$. Clearly, to get $f_{k-j}^{(2k)}$, $j =0,1,\ldots,k$, it suffices to consider the sub-matrix $M_j$ consisting of the first $j+1$ rows and the first $j+1$ columns of $M_k$. Thus, we may restrict ourselves to the following system of $j+1$ equations in the unknowns $f^{(2k)}_{k}, f^{(2k)}_{k-1}, \ldots, f^{(2k)}_{k-j}$,
\begin{equation}\label{sys1}
\setlength{\arraycolsep}{3.pt}
\begin{pmatrix}
\binom{2k+1}{1} & 0 & 0 & \hdots & 0 & 0 \\[6pt]
4^{-1}\binom{2k+1}{3} & \binom{2k-1}{1} & 0 & \hdots & 0 & 0 \\[8pt]
4^{-2}\binom{2k+1}{5} & 4^{-1}\binom{2k-1}{3} & \binom{2k-3}{1} & \hdots  & 0 & 0 \\[6pt]
\vdots & \vdots & \vdots & \ddots & \vdots & \vdots \\[7pt]
4^{-j+1}\binom{2k+1}{2j-1} & 4^{-j+2}\binom{2k-1}{2j-3} & 4^{-j+3}\binom{2k-3}{2j-5} & \hdots & \binom{2k-2j+3}{1} & 0 \\[7pt]
4^{-j}\binom{2k+1}{2j+1} & 4^{-j+1}\binom{2k-1}{2j-1} & 4^{-j+2}\binom{2k-3}{2j-3} & \hdots & 4^{-1}\binom{2k-2j+3}{3}
& \binom{2k-2j+1}{1} \\[5pt]
\end{pmatrix}
\begin{pmatrix} f_{k}^{(2k)} \\[6pt]  f_{k-1}^{(2k)} \\[6pt]  f_{k-2}^{(2k)} \\[2pt] \vdots \\[6pt]
  f_{k-j+1}^{(2k)} \\[6pt] f_{k-j}^{(2k)}
\end{pmatrix}
= \begin{pmatrix} 1\\[7pt] 0 \\[7pt] 0 \\[7pt] \vdots \\[7pt] 0 \\[7pt] 0
\end{pmatrix}.
\end{equation}
It is easily verified that the determinant of $M_j$ (which is equal to the product of the elements on the main diagonal) is given by
\begin{equation*}\label{det1}
|M_j | = \frac{(2k + 1)!!}{(2k-2j-1)!!},
\end{equation*}
where, for a positive integer $n$, the double factorial $n!!$ is defined as \cite{gould}
\begin{equation*}
n!! = (n)(n-2)(n-4) \cdots (4)(2),
\end{equation*}
if $n$ is even, and
\begin{equation*}
n!! = (n)(n-2)(n-4) \cdots (3)(1),
\end{equation*}
if $n$ is odd. In particular $|M_k| = (2k+1)!!$, where we assume the convention that $(-1)!! =1$. Since $|M_j| \neq 0$, Cramer's rule can be applied to the system \eqref{sys1} to obtain the following determinant formula for $f^{(2k)}_{k-j}$
\begin{equation}\label{det2}
f^{(2k)}_{k-j} = (-1)^j \frac{(2k-2j-1)!!}{(2k+1)!!} \Delta_j, \quad j=0,1,\ldots,k,
\end{equation}
where $\Delta_0 =1$ and, for $j =1,\ldots,k$, $\Delta_j$ is the determinant of order $j$
\begin{equation}\label{det3}
\Delta_j =
\setlength{\arraycolsep}{3.pt}
\begin{vmatrix}
4^{-1}\binom{2k+1}{3} & \binom{2k-1}{1} & 0 & \hdots & 0  \\[8pt]
4^{-2}\binom{2k+1}{5} & 4^{-1}\binom{2k-1}{3} & \binom{2k-3}{1} & \hdots  & 0 \\[6pt]
\vdots & \vdots & \vdots & \ddots & \vdots \\[7pt]
4^{-j+1}\binom{2k+1}{2j-1} & 4^{-j+2}\binom{2k-1}{2j-3} & 4^{-j+3}\binom{2k-3}{2j-5} & \hdots & \binom{2k-2j+3}{1} \\[7pt]
4^{-j}\binom{2k+1}{2j+1} & 4^{-j+1}\binom{2k-1}{2j-1} & 4^{-j+2}\binom{2k-3}{2j-3} & \hdots & 4^{-1}\binom{2k-2j+3}{3} \\[5pt]
\end{vmatrix}.
\end{equation}

Similarly, starting from \eqref{odd} and the relation
\begin{equation*}
S_{2k+1}(n) - S_{2k+1}(n-1) = n^{2k+1},
\end{equation*}
we arrive at the following linear triangular system of $k+1$ equations in the unknowns $f^{(2k+1)}_{k}, f^{(2k+1)}_{k-1}, \ldots, f^{(2k+1)}_{0}$
\begin{equation}\label{recur2}
\sum_{m = j}^{k} 4^{j-m} \binom{2m+2}{2j+1} f^{(2k+1)}_{m} = \delta_{j,k}, \quad j =0,1,\ldots,k.
\end{equation}
Hence, proceeding as above, we can derive the corresponding determinant formula for $f^{(2k+1)}_{k-j}$, namely
\begin{equation}\label{det4}
f^{(2k+1)}_{k-j} = (-1)^j \frac{(2k-2j)!!}{(2k+2)!!} \Delta_j^{\prime}, \quad j=0,1,\ldots,k,
\end{equation}
where $\Delta_0^{\prime} =1$ and, for $j =1,\ldots,k$, $\Delta_j^{\prime}$ is the determinant of order $j$
\begin{equation}\label{det5}
\Delta_j^{\prime} =
\setlength{\arraycolsep}{3.pt}
\begin{vmatrix}
4^{-1}\binom{2k+2}{3} & \binom{2k}{1} & 0 & \hdots & 0  \\[8pt]
4^{-2}\binom{2k+2}{5} & 4^{-1}\binom{2k}{3} & \binom{2k-2}{1} & \hdots  & 0 \\[6pt]
\vdots & \vdots & \vdots & \ddots & \vdots \\[7pt]
4^{-j+1}\binom{2k+2}{2j-1} & 4^{-j+2}\binom{2k}{2j-3} & 4^{-j+3}\binom{2k-2}{2j-5} & \hdots & \binom{2k-2j+4}{1} \\[7pt]
4^{-j}\binom{2k+2}{2j+1} & 4^{-j+1}\binom{2k}{2j-1} & 4^{-j+2}\binom{2k-2}{2j-3} & \hdots & 4^{-1}\binom{2k-2j+4}{3} \\[5pt]
\end{vmatrix}.
\end{equation}

\subsection{Example}

For $k=5$, from \eqref{det2} and \eqref{det3} we obtain
\begin{align*}
f_5^{(10)} & = \frac{1}{11}  \\
f_4^{(10)} & = - \frac{7!!}{11!!} 4^{-1} \binom{11}{3} = - \frac{5}{12} \\
f_3^{(10)} & = \frac{5!!}{11!!}
\begin{vmatrix}
4^{-1}\binom{11}{3} & \binom{9}{1} \\[6pt]
4^{-2}\binom{11}{5} & 4^{-1}\binom{9}{3}\\
\end{vmatrix} = \frac{7}{8} \\
f_2^{(10)} & = - \frac{3}{11!!}
\begin{vmatrix}
4^{-1}\binom{11}{3} & \binom{9}{1} & 0 \\[6pt]
4^{-2}\binom{11}{5} & 4^{-1}\binom{9}{3} & \binom{7}{1} \\[6pt]
4^{-3}\binom{11}{7} & 4^{-2}\binom{9}{5} & 4^{-1}\binom{7}{3}
\end{vmatrix} = -\frac{31}{32} \\
f_1^{(10)} & = \frac{1}{11!!}
\begin{vmatrix}
4^{-1}\binom{11}{3} & \binom{9}{1} & 0 & 0 \\[6pt]
4^{-2}\binom{11}{5} & 4^{-1}\binom{9}{3} & \binom{7}{1} & 0 \\[6pt]
4^{-3}\binom{11}{7} & 4^{-2}\binom{9}{5} & 4^{-1}\binom{7}{3} & \binom{5}{1} \\[6pt]
4^{-4}\binom{11}{9} & 4^{-3}\binom{9}{7} & 4^{-2}\binom{7}{5} & 4^{-1}\binom{5}{3}
\end{vmatrix} = \frac{127}{256}  \\
f_0^{(10)} & = -\frac{1}{11!!}
\begin{vmatrix}
4^{-1}\binom{11}{3} & \binom{9}{1} & 0 & 0 & 0 \\[6pt]
4^{-2}\binom{11}{5} & 4^{-1}\binom{9}{3} & \binom{7}{1} & 0 & 0 \\[6pt]
4^{-3}\binom{11}{7} & 4^{-2}\binom{9}{5} & 4^{-1}\binom{7}{3} & \binom{5}{1} & 0 \\[6pt]
4^{-4}\binom{11}{9} & 4^{-3}\binom{9}{7} & 4^{-2}\binom{7}{5} & 4^{-1}\binom{5}{3} & \binom{3}{1} \\[6pt]
4^{-5}\binom{11}{11} & 4^{-4}\binom{9}{9} & 4^{-3}\binom{7}{7} & 4^{-2}\binom{5}{5} & 4^{-1}\binom{3}{3}
\end{vmatrix} = -\frac{2555}{33792}
\end{align*}
where the determinants have been computed using \emph{Mathematica} software. Then we have
\begin{equation*}
S_{10}(n) = \frac{1}{11}N^{11} - \frac{5}{12}N^9 + \frac{7}{8} N^7 - \frac{31}{32} N^5 + \frac{127}{256} N^3
- \frac{2555}{33792}N,
\end{equation*}
where $N$ is a shorthand for $n + \frac{1}{2}$.

On the other hand, to get the Faulhaber polynomial $S_{11}(n)$, we can of course use the equations \eqref{det4} and \eqref{det5} to find $f_m^{(11)}$, $m=0,1,\ldots,5$, and then using \eqref{ind} to find $c_{11}$. However, there is a handy shortcut for finding $f_m^{(11)}$ \emph{once} the corresponding coefficient $f_m^{(10)}$ is known. To see this in full generality, we invoke the following elementary result according to which \cite{sherwood,owens,wu}
\begin{equation*}
S_{2k+1}^{\prime}(n) = \frac{\text{d}S_{2k+1}(n)}{\text{d}n} = (2k+1) S_{2k}(n), \quad k \geq 1,
\end{equation*}
where, momentarily, we treat $n$ as a continuous variable for convenience. Thus, taking the derivative of \eqref{odd} with respect to $n$ yields
\begin{equation*}
S_{2k+1}^{\prime}(n) = \sum_{m=0}^{k} (2m+2) f^{(2k+1)}_{m} \left( n+\frac{1}{2} \right)^{2m+1},
\end{equation*}
and then
\begin{equation*}
S_{2k}(n) = \sum_{m=0}^{k} \frac{2m+2}{2k+1} f^{(2k+1)}_{m} \left( n+\frac{1}{2} \right)^{2m+1}.
\end{equation*}
Comparing the last equation with \eqref{even}, and noting that $\{ (n+ 1/2)^{2m+1} \}_{m=0}^{k}$ constitutes a set of $k+1$ linearly independent polynomials, we conclude that
\begin{equation}\label{iso}
f_m^{(2k+1)} = \frac{2k+1}{2m+2} f_m^{(2k)}, \quad m =0,1,\ldots,k \,\,\, \text{and}\,\,\,  k\geq 1.
\end{equation}
Now, employing \eqref{iso} and the values of $f_0^{(10)}, \ldots, f_5^{(10)}$ previously obtained, we get
\begin{equation*}
S_{11}(n) = \frac{1}{12}N^{12} - \frac{11}{24}N^{10} + \frac{77}{64} N^8 - \frac{341}{192} N^6 + \frac{1397}{1024} N^4
- \frac{2555}{6144}N^2 + \frac{691}{16384},
\end{equation*}
where the last term $c_{11} = \frac{691}{16384}$ is obtained from \eqref{ind}.

\subsection{Witmer's recursive formulas}

In a largely unnoticed paper, \citeA[Theorems III and IV]{witmer} derived recursive formulas for the coefficients $f^{(2k)}_{m}$, $f^{(2k+1)}_{m}$, and $c_{2k+1}$. In our notation, these formulas read as
\begin{align*}
f_{k}^{(2k)} & = \frac{1}{2k+1},  \\
f_{i}^{(2k)} & = -\frac{1}{2k+1} \sum_{j=1}^{k-1} 4^{j-k} \binom{2k+1}{2j} f_{i}^{(2j)}, \quad 1 \leq i < k  \\
f_{0}^{(2k)} & = -\frac{1}{2k+1}\left[ \frac{1}{4^k} + \sum_{j=1}^{k-1} 4^{j-k} \binom{2k+1}{2j} f_{0}^{(2j)} \right],
\end{align*}
and
\begin{align*}
f_{k}^{(2k+1)} & = \frac{1}{2k+2},  \\
f_{i}^{(2k+1)} & = -\frac{1}{2k+2} \sum_{j=0}^{k-1} 4^{j-k} \binom{2k+2}{2j+1} f_{i}^{(2j+1)}, \quad 0 \leq i < k  \\
c_{2k+1} & = -\frac{1}{2k+2}\left[ \frac{1}{4^{k+1}} + \sum_{j=0}^{k-1} 4^{j-k} \binom{2k+2}{2j+1} c_{2j+1} \right].
\end{align*}
It is to be noted that, in contrast to our \emph{horizontal} recursive formulas in equations \eqref{recur1} and \eqref{recur2}, Witmer's \emph{vertical} recursive formula for $f_{i}^{(2k)}$ [$f_{i}^{(2k+1)}$] does require the knowledge of the Faulhaber coefficients $f_{i}^{(2j)}$ [$f_{i}^{(2j+1)}$] corresponding to the power sum polynomials $S_{2j}(n)$ [$S_{2j+1}(n)$], for $j =1,\ldots, k-1$ [$j=0,1,\ldots, k-1$]. Likewise, Witmer's formula for $c_{2k+1}$ requires knowing the previous coefficients $c_1, c_3, \ldots, c_{2k-1}$, for $k \geq 1$.

\subsection{Explicit formulas}

For completeness, next we write down an explicit formula for the coefficients $f^{(2k)}_{m}$ and $f^{(2k+1)}_{m}$ in terms of Bernoulli numbers, namely \cite{cere2}
\begin{align}
f^{(2k)}_{m} & = \frac{1}{2m+1} \binom{2k}{2m} B_{2k-2m}\left( \frac{1}{2}\right), \label{expl1} \\
\intertext{and}
f^{(2k+1)}_{m} & = \frac{1}{2m+2} \binom{2k+1}{2m+1} B_{2k-2m}\left( \frac{1}{2}\right),  \label{expl2}
\end{align}
for $m =0,1,\ldots,k$, and where $B_r \left( \frac{1}{2}\right) = \big( 2^{1-r} -1 \big)B_r$ is the value taken by the Bernoulli polynomial $B_r(x)$ at $ x = \frac{1}{2}$. In particular, $f_0^{(2k)} = B_{2k}\left( \frac{1}{2}\right)$ and $f_0^{(2k+1)} = \frac{2k+1}{2}B_{2k}\left( \frac{1}{2}\right)$. Thus, setting $j =k$ in \eqref{det2} and \eqref{det4}, and solving respectively for $\Delta_k$ and $\Delta_k^{\prime}$, gives
\begin{equation*}
\Delta_k = (-1)^k (2k+1)!! B_{2k}\left( \frac{1}{2}\right),
\end{equation*}
and
\begin{equation*}
\Delta_k^{\prime} = (-1)^k (2k+1)(k+1)(2k)!! B_{2k}\left( \frac{1}{2}\right),
\end{equation*}
for $k =1,2,3,\ldots\,$.

Moreover, using \eqref{expl1} and \eqref{expl2}, it is not difficult to show the following explicit formulas for the Faulhaber polynomials \eqref{even} and \eqref{odd}, namely
\begin{equation*}
S_{2k}(n) = \frac{1}{2k+1} \sum_{j=0}^{k} \binom{2k+1}{2j} B_{2j} \left( \frac{1}{2}\right) N^{2k+1-2j}, \quad{k \geq 1},
\end{equation*}
and
\begin{equation*}
S_{2k+1}(n) = \frac{1}{2k+2} \sum_{j=0}^{k} \binom{2k+2}{2j} B_{2j} \left( \frac{1}{2}\right) \left[ N^{2k+2-2j}
- 4^{j-k-1} \right],  \quad{k \geq 0}.
\end{equation*}

\section{Concluding remarks}

As explained in \citeA{hersh}, the form of the Faulhaber polynomials \eqref{even} and \eqref{odd} stems from the fact that $S_p(n)$ is symmetric about $-\frac{1}{2}$. More precisely, $S_p(n)$ fulfills the symmetry property \cite{levy,krishna,shirali,newsome}
\begin{equation}\label{symm}
S_p ( -(n+1)) = (-1)^{p+1} S_p(n),  \quad p \geq 1.
\end{equation}
Therefore, for even $p =2k$, $k \geq 1$, we have $S_{2k}(n) = - S_{2k}(-(n+1))$, which implies that $S_{2k}(n)$ is symmetric about the point $(-\frac{1}{2},0)$. Hence, $S_{2k}(n)$ can be expressed as an odd polynomial in $n + \frac{1}{2}$, in accordance with \eqref{even}. On the other hand, for odd $p=2k+1$, $k \geq 0$, we have $S_{2k+1}(n) = S_{2k+1}(-(n+1))$, which implies that $S_{2k+1}(n)$ is symmetric about the vertical line at $-\frac{1}{2}$. This in turn means that $S_{2k+1}(n)$ can be expressed as a even polynomial in $n + \frac{1}{2}$, in accordance with \eqref{odd}. Note also that $S_{2k+1}(-\frac{1}{2}) = c_{2k+1}$. It is to be mentioned, on the other hand, that \eqref{symm} is equivalent to the symmetry property of the Bernoulli polynomials \cite{apostol}
\begin{equation*}
B_p (1-x) = (-1)^{p} B_p(x),  \quad p \geq 0.
\end{equation*}
This is so due to: (i) the relationship between $S_p(n)$ and the Bernoulli polynomials, namely, $S_p(n) = \frac{B_{p+1}(n+1) - B_{p+1}}{p+1}$, $p \geq 1$; and (ii) the fact that $B_{2p+1}=0$ for all $p \geq 1$.

To conclude, it should be emphasized that the Faulhaber polynomials \eqref{even} and \eqref{odd} can be expressed equivalently in the form \cite{edwards1,edwards2,knuth,beardon,dubeau,cere}
\begin{align}
S_{2k}(n) & = S_2(n) \big[ b_{k,0} + b_{k,1} S_1(n) + b_{k,2} \big( S_1(n) \big)^2 + \cdots +
b_{k,k-1} \big( S_1(n) \big)^{k-1} \big], \label{even2} \\[-2mm]
\intertext{and}
S_{2k+1}(n) & =  \big( S_1(n)\big)^2 \big[ c_{k,0} + c_{k,1} S_1(n) + c_{k,2} \big( S_1(n) \big)^2 + \cdots +
c_{k,k-1} \big( S_1(n) \big)^{k-1} \big], \label{odd2}
\end{align}
respectively, where $b_{k,j}$ and $c_{k,j}$ are non-zero rational coefficients for $j =0,1,\ldots,\linebreak k-1$ and $k \geq 1$. Indeed, noting that $\big( n + \frac{1}{2} \big)^2 = \frac{1}{4} \big( 1 + 8 S_1(n) \big)$, and identifying the polynomials \eqref{even} with \eqref{even2}, and \eqref{odd} with \eqref{odd2}, it can be shown that
\begin{align}
b_{k,j} & = \frac{3}{2} \, 8^{j+1} \sum_{m= j+1}^{k} \binom{m}{j+1} \frac{f_m^{(2k)}}{4^m}, \label{id1} \\[-2mm]
\intertext{and}
c_{k,j} & = 8^{j+2} \sum_{m= j+1}^{k} \binom{m+1}{j+2} \frac{f_m^{(2k+1)}}{4^{m+1}}. \label{id2}
\end{align}
Notice that \eqref{id1} [\eqref{id2}] does not involve the coefficient $f_0^{(2k)}$ [$f_0^{(2k+1)}$]. Thus, plugging the known values of $f_m^{(10)}$ and $f_m^{(11)}$, $m =1,\ldots,5$, into \eqref{id1} and \eqref{id2}, respectively, leads to
\begin{align*}
S_{10}(n) & = S_2(n) \left[ \frac{5}{11} - \frac{30}{11} S_1(n) + \frac{68}{11} \big( S_1(n) \big)^2 - \frac{80}{11}
\big( S_1(n) \big)^3 + \frac{48}{11} \big( S_1(n) \big)^4 \right],  \\
\intertext{and}
S_{11}(n) & = \big( S_1(n)\big)^2  \left[ \frac{5}{3} - \frac{20}{3} S_1(n) + \frac{34}{3} \big( S_1(n) \big)^2 - \frac{32}{3}
\big( S_1(n) \big)^3 + \frac{16}{3} \big( S_1(n) \big)^4 \right].
\end{align*}
Furthermore, recalling \eqref{iso}, from \eqref{id1} and \eqref{id2} it follows that the coefficients $b_{k,j}$ and $c_{k,j}$ are related by
\begin{equation*}
c_{k,j} = \frac{4k+2}{3j+6} \, b_{k,j},   \quad j =0,1,\ldots,k-1 \,\,\,\text{and}\,\,\, k\geq 1.
\end{equation*}

Conversely, one can get the set of coefficients $\{ f^{(2k)}_{m}\}_{m=0}^{k}$ and $\{f^{(2k+1)}_{m} \}_{m=0}^{k}$ from $\{ b_{k,j} \}_{j=0}^{k-1}$ and $\{ c_{k,j} \}_{j=0}^{k-1}$, respectively, by means of the relations
\begin{align*}
f_m^{(2k)} & = \frac{2}{3} \, (-4)^m  \sum_{j= m}^{k} \frac{(-1)^{j}}{8^{j}} \binom{j}{m} b_{k,j-1}, \\[-2mm]
\intertext{and}
f_m^{(2k+1)} & = (-4)^{m+1}  \sum_{j= m}^{k} \frac{(-1)^{j+1}}{8^{j+1}} \binom{j+1}{m+1} c_{k,j-1},
\end{align*}
where it is understood that $b_{k,-1} = c_{k, -1} =0$. On the other hand, the term $c_{2k+1}$ is obtained as the value of the polynomial \eqref{odd2} at $n = - \frac{1}{2}$, namely
\begin{equation*}
c_{2k+1} = \sum_{j=0}^{k-1} \frac{(-1)^j}{8^{j+2}} \, c_{k,j}.
\end{equation*}

\bibliography{references2}
\bibliographystyle{apacite}
	
\end{document}